\documentclass[12pt, letterpaper,reqno]{amsart}

\title{Topological K-Theory of Complex Projective Spaces}
\author{Virgil Chan}
\date{28 February 2013}
\email{virchan@ucdavis.edu}

\usepackage{tikz,amsmath,setspace,amssymb,enumerate,subfig,indentfirst,fullpage,hyperref,framed,etoolbox,bm,fancyhdr,lipsum,footnote}

\makeatletter
\patchcmd{\@maketitle}
  {\ifx\@empty\@dedicatory}
  {\ifx\@empty\@date \else {\vskip3ex \centering\footnotesize\@date\par\vskip1ex}\fi
   \ifx\@empty\@dedicatory}
  {}{}
\patchcmd{\@adminfootnotes}
  {\ifx\@empty\@date\else \@footnotetext{\@setdate}\fi}
  {}{}{}
\makeatother

\newtheoremstyle{break}
  {\topsep}{\topsep}%
  {\itshape}{}%
  {\bfseries}{}%
  {\newline}{}%

\theoremstyle{break}
\newtheorem{theorem}{Theorem}[section]
\newtheorem{lemma}[theorem]{Lemma}
\newtheorem{proposition}[theorem]{Proposition}
\newtheorem{corollary}[theorem]{Corollary}
\newtheorem{definition}[theorem]{Definition}

\newcommand{\D}[1]{\mathbb#1}
\newcommand{\sphere}{\D{S}}
\newcommand{\PC}{\D{C}\mbox{P}}
\newcommand{\Ext}{\mbox{Ext}}
\newcommand{\CG}{\D{C}\mbox{G}}

\usetikzlibrary{arrows,chains,matrix,positioning,scopes}

\makeatletter
\tikzset{join/.code=\tikzset{after node path={%
\ifx\tikzchainprevious\pgfutil@empty\else(\tikzchainprevious)%
edge[every join]#1(\tikzchaincurrent)\fi}}}
\makeatother

\tikzset{>=stealth',every on chain/.append style={join},
         every join/.style={->}}

\newlength{\parindentsave}

\numberwithin{equation}{section}

\begin{document}

\begin{abstract}
We compute the $K$-theory of complex projective spaces. There are three major ingredients: the exact sequence of $K$-groups, the theory of Chern character and the Bott Periodicity Theorem.
\end{abstract}

\maketitle

\tableofcontents

\section{Introduction}

Topological $K$-theory (or $K$-theory in short) is the study of abelian groups generated by vector bundles. It is an extraordinary cohomology theory that plays an important role in topology. The fundamental concept of $K$-theory is the construction of the Grothendieck group (see Proposition \ref{grothendieck group}) from the equivalence classes of complex vector bundles. That is, for any complex vector bundle, we associate it with a sequence of abelian groups, known as the $K$-groups or $K$-functor. According to the common opinion, it was Alexander Grothendieck who had started the subject to formulate his Grothendieck-Riemann-Roch Theorem, but the first works in $K$-theory were published in 1959 by Michael Atiyah and Friedrich Hizebruch, and in 1964 most results were completed due to Frank Adams, Michael Atiyah, Raoul Bott and Freidrich Hizebruch. The most remarkable result related to $K$-theory perhaps is the Atiyah-Singer Index Theorem in Partial Differential Equations. Later, the key ideas of $K$-theory was extended to algebra and algebraic geometry, which nowadays known as the Algebraic $K$-theory (see \cite{rosenberg}).

The goal of this paper is to compute the $K$-theory of the complex projective space $\PC^n$:
\newsavebox{\maintheorem}\begin{lrbox}{\maintheorem}
\begin{minipage}{\textwidth}\setlength{\parindent}{\parindentsave}
\vspace{\abovedisplayskip}
\begin{leftbar}
\begin{theorem}
\label{hahaj}
\[ K^q(\PC^n) = \left\{ \begin{array}{ll}
         \D{Z}^{\oplus n+1} & \mbox{if $q$ is even};\\
        0 & \mbox{if $q$ is odd}.\end{array} \right. \]
        
Moreover, as a ring, 
 \[K^0(\PC^n) = \D{Z}[\gamma]/\langle \gamma^{n+1} \rangle, \] where $\gamma = \zeta - \tilde{1}$, $\zeta$ is the Hopf bundle, and $\tilde{1}$ is the trivial vector bundle $\PC \times \D{C} \rightarrow \PC$ of dimension $1$.       
\end{theorem}
\end{leftbar}
\end{minipage}\end{lrbox}\newline\noindent\usebox{\maintheorem}\\

In Section 2, we introduce the basic notations and concepts of $K$-theory based on \cite{fuchs}.

In Section 3, we study the cohomology structure of $\PC^n$ by using its CW structure.

In Section 4, we briefly introduce the theory of Chern character based on \cite{milnorsta} and \cite{madsen}, and compute the Chern character of the Hopf bundle.

In Section 5, we state the Bott Periodicity Theorem and its corollaries to path the basic settings for the proof of Theorem \ref{hahaj}.

\subsection*{Acknowledgement} The author is extraordinarily grateful to his undergraduate thesis advisor \href{http://www.math.ucdavis.edu/research/profiles/fuchs}{Professor Dmitry Fuchs} for his uncountable help, insightful discussions and endless encouragement in writing this paper. Nevertheless, the author would also like to thank \href{http://www.math.ucdavis.edu/research/profiles/schwarz}{Professor Albert Schwarz} and \href{http://www.math.ucdavis.edu/~kaminker/}{Professor Jerome Kaminker} in reviewing this paper and providing fruitful feedback.

\section{K-Theory}

We start with a proposition:
\begin{leftbar}
\begin{proposition}[{{\bf Construction of the Grothendieck Group}}]
\label{grothendieck group}
Let $S$ be a commutative semi-group (not necessarily having a unit). There is an abelian group $G$ (called the \textbf{Grothendieck group} or \textbf{group completion} of $S$) and a homomorphism $\phi: S \rightarrow H$, such that for any group abelian group $H$ and homomorphism $\psi:S\to H$, there is a unique homomorphism $\theta:G\to H$ such that $\psi=\theta \circ \phi$.
\end{proposition}
\end{leftbar}

\begin{proof}
Define an equivalence relation ``$\sim$'' in $S \times S$ by:
\[
(x,y)\sim(u,v) \iff \exists \ t\in S \ \mbox{such that } x+v+t=u+y+t \ \mbox{in} \ S.
\]

Denote the equivalence class of $(x,y)$ by $\left[ (x,y) \right]$, and let $G$ be the set of equivalence classes of ``$\sim$''. We define the abelian group structure in $G$.

First of all, the addition of equivalence classes is defined to be:
\[ \left[ (x,y) \right] + \left[ (x',y') \right] = \left[ (x+x', y+y') \right], \]
and we show that it is well-defined.

If $(x,y) \sim (u,v)$ and $(x',y') \sim (u',v')$, then by definition, there exists $t$, $t' \in S$ such that:
\begin{eqnarray*}
x+v+t &=& y+u+t; \\
x'+v'+t' &=& y'+u'+t'.
\end{eqnarray*}
So we have:
\begin{eqnarray*}
(x+x')+(v+v')+(t+t') &=& (x+v+t)+(x'+v'+t') \\
&=& (y+u+t)+(y'+u'+t') \\
&=& (y+y') +(u+u')+(t+t')
\end{eqnarray*}
and hence $\left[ (x+x',y+y') \right] \sim \left[ (u+u'), (v+v') \right]$. This concludes the addition is well-defined.

Next, we define the identity element with respect to this addition.

Note that for any $x$, $y \in S$, we have:
\[ \left[ (x,x) \right] = \left[ (y,y) \right] \]
since $x+y = y+x$. We denote this distinguished element $\left[ (x,x) \right]$ by $0$. Moreover, $0$ is indeed the additive identity since for every $u$, $v \in S$, we have:
\[ (u+x, v+x) \sim (u,v). \]

Equally important, we wish to define the additive inverse.

Just observe that for every $x$, $y \in S$, we have:
\[ \left[ (x,y) \right] + \left[ (y,x) \right] = \left[ (x+y, x+y) \right] = 0.\]
This shows that $\left[ (x,y) \right] = -\left[ (y,x) \right]$

If $(x,y) \sim (u,v)$, then there exists $t$ such that:
\[ x+v+t = y+u+t. \]
$\Rightarrow y+u+t = x+v+t$ and hence, $(y,x) \sim (v,u)$.

It is easy to check also that the addition, together with the inverse and identity satisfy the group axioms. Moreover, commutativity of elements in $G$ follows from the the commutative property of the semi-group $S$. As a result, $G$ is an abelian group as desired.

Now, define $\phi: S\to G$ such that $\phi (x)=[(x+x,x)]$. Note that $[(x,y)]=\phi (x)- \phi (y)$. Therefore Im($\phi$) generates $G$. 

Given a group $H$ and homomorphism $\psi: S\to H$, define $\theta:G\to H$ by $\theta([(x,y)])= \psi(x)-\psi(y)$. We have $\psi=\theta \circ \phi$.

If $\phi':S\to G'$ is any other pair with same property, then there exists an isomorphism $\alpha:G\to G'$ such that $\phi'= \alpha \circ \phi$.

Therefore, the claim holds as desired.
\end{proof}

Now, consider a finite CW complex $X$. Let $F(X)$ be the set of equivalence classes of complex vector bundles with the base $X$. Observe that $F(X)$ has two binary operations: the Whitney sum ``$\oplus$'' and the tensor product ``$\otimes$'', correspond to addition and multiplication respectively. Also note that the class of zero-dimensional complex vector bundle serves as the additive identity element in $F(X)$. So $F(X)$ forms a semi-group under $\oplus$. By Proposition \ref{grothendieck group}, we can construct an abelian group from $F(X)$. This motivates the following definition:

\begin{leftbar}
\begin{definition}[{{\bf $\bm{K}$-group}}]
\label{def of K-group}
Let $X$ be a finite CW complex. Let $F(X)$ be the set of equivalence classes of complex vector bundles with the base $X$. Then, $K(X)$ is the Grothendieck group of $F(X)$. We call elements in $K(X)$ as {\bf virtual bundles}.
\end{definition}
\end{leftbar}

We describe a little bit more of the algebraic structure of $K(X)$. First, direct sum determines the addition in $K(X)$, with the equivalence class of zero-dimensional vector bundle as the additive identity. Secondly, the tensor product determines a multiplication in $K(X)$, with the equivalence class of trivial bundle of dimension $1$ is the multiplicative identity in $K(X)$. Moreover, one can easily show that the multiplication in $K(X)$ is commutative, associative and distributive with respect to the addition. Therefore, we can actually answer questions concerning the ring structure of $K(X)$.

Consider a trivial case when $X$ is a single point. Clearly, a complex vector bundle over a point is uniquely determined by its dimension. So the dimension map, denote ``dim", gives an isomorphism from $F(\mbox{point})$ to $\D{N}$. That is, $F(\mbox{point}) = \D{N}$. Since the Grothendieck group of $\D{N}$ is $\D{Z}$, we have:
\[ K(\mbox{point}) = \D{Z}. \]

Similarly, for any finite CW complex $X$, the dimension map gives a map
\[ \mbox{dim}:F(X) \rightarrow \D{N} \]
by sending the trivial bundle $\tilde{1}$ to the number $1 \in \D{Z}$, which extended to a group homomorphism:
\[ \mbox{dim}:K(X) \rightarrow \D{Z} \]
and leads to the following definition:

\begin{leftbar}
\begin{definition}[{{\bf Reduced ${\bm K}$-group ${\bm \tilde{{\bm K}}}$}}]
\label{def of reduced K-group}
Let $X$ be a finite CW complex. The {\bf reduced ${\bm K}$-group of ${\bm X}$} is defined to be:
\[ \tilde{K}(X) = \mbox{Ker} \left[ \mbox{dim}:K(X) \rightarrow \D{Z} \right] \]
where ``dim" is the dimension map.

Equivalently, we have:
\[ \tilde{K}(X) = \mbox{Ker}\left[ \mbox{dim}:K(X) \rightarrow K(x_0) \right], \]
where $x_0$ is the base-point of $X$.
\end{definition}
\end{leftbar}

A geometrical way to understand $\tilde{K}$ is that, $\tilde{K}(X)$ is the group of stably equivalent classes of complex vector bundles over $X$ (see \cite[Chapter 6, Section 37, P.4 Theorem]{fuchs}).

\begin{leftbar}
\begin{definition}[{{\bf Relative ${\bm K}$-group ${\bm K^{{\bm q}}}\bm{(X,A)}$}}]
\label{def of higher K-groups and relative K-groups}
Let $(X,A)$ be a finite CW pair and $\Sigma X$ denote the suspension over $X$. For $q \leq 0$, we define:
\[ K^q(X,A) = \tilde{K}(\Sigma^{-q}(X/A)). \]
\end{definition}
\end{leftbar}

From the above definition, if we pick $A = \emptyset$, then:
\[ X/A = X/\emptyset = X \sqcup (\mbox{one point}). \]

Hence,
\[ K(X,\emptyset) = \tilde{K}(X/\emptyset) = \tilde{K}(X \sqcup \ \mbox{point}) = K(X) \]
because $\tilde{K}(X \sqcup \mbox{point})$ is the group of virtual bundles over $X\sqcup\mbox{point}$, that is dimension $0$ over the point, which is exactly the group of virtual bundles over $X$.
Therefore, we have:
\begin{equation}
K^q(X) = K^q(X,\emptyset).
\end{equation}

Recall that if $\xi$ is a vector bundle $p: \mathcal{E} \rightarrow Y$, and~$f:X \rightarrow Y$ is a continuous map, then the {\bf induced vector bundle} is defined to be:
\[ f^{\ast}\xi = \lbrace (x,e) \in X \times \mathcal{E} \ \vert p(e) = f(x) \rbrace. \]
Moreover, the vector bundle $f^{\ast}\xi$ over $X$ has same dimension of the vector bundle $\mathcal{E}$ over $Y$.

The construction of the induced bundle establishes a ring homomorphism $K(Y) \rightarrow K(X)$ (and hence also $\tilde{K}(Y) \rightarrow \tilde{K}(X)$ and $K^{-q}(Y) \rightarrow K^{-q}(X)$). With this homomorphism, one can show that the functor $\tilde{K}$ is half exact (see \cite{fuchs}), and for a CW pair $(X,A)$, we have the following {\bf exact sequence of $K$-groups}:

\begin{equation}
\begin{tikzpicture}[node distance=6mm]
\label{exact sequence of k groups}
  \node (A) {$K^0(A)$}; 
  \node (B) [right= of A] {$\cdots$}; 
  \node (C) [right= of B] {$K^{-q}(A)$}; 
  \node (D) [right= of C] {$K^{-q}(X)$};
  \node (E) [right= of D] {$K^{-q}(X,A)$};
  \node (F) [right= of E] {$K^{-q-1}(A)$};
  \node (G) [right= of F] {$\cdots$};
  
  \draw[stealth-] (A)-- node[right] {\small} (B); 
  \draw[stealth-] (B)-- node [right] {\small} (C); 
  \draw[stealth-] (C)-- node [right] {\small} (D);
  \draw[stealth-] (D)-- node [right] {\small} (E);
  \draw[stealth-] (E)-- node [right] {\small} (F);
  \draw[stealth-] (F)-- node [right] {\small} (G);
\end{tikzpicture}
\end{equation}
which will be the one of the main tools to prove Theorem \ref{hahaj} and complete the construction of $K$-groups.

\subsection{K-theory and Classifying Space}\mbox{}\\

The general linear group $GL(R)$ (where $R$ is a ring) arises in lots of areas in mathematics. In the context of topological and algebraic $K$-theory, it is an essential element. For example, in algebraic $K$-theory, the definition of the first algebraic $K$-group $K_1$ is defined to be the quotient group of $GL(R)$ and its subgroup generated by elementary matrices (see \cite{rosenberg}). In topological $K$-theory, its subgroup $U$ allows one to compute the $K$-theory of certain CW complexes.

\begin{leftbar}
\begin{definition}[{\bf The group $\bm{U}$}]
Let $U(n)$ be the group of unitary $n\times n$-matrices. We embed $U(n)$ into $U(n+1)$ by:
\[ A\mapsto \left[ \begin{array}{ccc}
A & 0 \\
0 & 1
\end{array} \right].\] 

The group $U$ is defined to be:
\[ U=\displaystyle \lim_{n \rightarrow \infty} U(n). \]
\end{definition}
\end{leftbar}

\begin{leftbar}
\begin{definition}[{\bf Classifying Space}]
Let $\CG(N,n)$ denote the complex Grassmannian. We can embed $\CG(N,n)$ to a larger Grassmannian $\CG(N',n')$ by a similar way we did to $U(n)$. The {\bf Classifying Space} is defined to be:
\[ BU = \CG(\infty, \infty) = \displaystyle \lim_{n \rightarrow \infty} \lim_{N \rightarrow \infty} \CG(N,n),\]
\end{definition}
\end{leftbar}

We now state an important relation between $K$-theory and Classifying Space.

\begin{leftbar}
\begin{theorem}[{{\cite[{\bf Chapter 6, Section 37, Corollary}]{fuchs}}}]
\label{k theory and bu}
For any finite CW complex $X$, $\tilde{K}(X)$ is equal to the set of homotopic maps from $X$ to $BU$. That is,
\[ \tilde{K}(X) = \pi(X,BU). \]
Moreover, for a continuous map $f:X \longrightarrow Y$, the induced homomorphism:
\[ f^{\ast}:\tilde{K}(Y) \longrightarrow \tilde{K}(X)\]
coincides with the map: 
\[ f^{\#}: \pi(Y,BU) \longrightarrow \pi(X,BU)\]
\end{theorem}
\end{leftbar}

\begin{proof}
As stated before, $\tilde{K}(X)$ is the group of stably equivalent classes of complex vector bundles over $X$. We show that the classes are one-to-one correspond to homotopy classes of the maps $X \rightarrow BU$.

By elementary bundle theory, a vector bundle of dimension $n$ over $X$ induces a map: 
\[ X \rightarrow \CG(\infty, n). \]
Since $\CG(\infty, n) \subset \CG(\infty, \infty)$, we get a continuous map $X \rightarrow \CG(\infty, \infty)$.

It follows immediately that two maps that correspond to two vector bundles over $X$
\[ X \rightarrow \CG(\infty, n_1), X \rightarrow \CG(\infty, n_2) \]
are homotopic in $\CG(\infty, \infty)$ if and only if the two vector bundles are stably equivalent. So the claim holds as desired.
\end{proof}

Using Theorem \ref{k theory and bu}, we are able to compute the $K$-theory of sphere. Put $X = \sphere^r$, we immediately get:
\[ \tilde{K}(\sphere^r) = \pi(\sphere^r, BU) = \pi_r(BU). \]
We can say more than this. By elementary homotopy theory, $U$ is homotopic equivalent to the loop space of $BU$ (i.e. $U \sim \Omega BU$). So a $(r-1)$-spheroid in $U$ is a $r$-spheroid in $BU$. It follows that $\pi_r(BU) = \pi_{r-1}(U)$. Gluing all these together we get:

\begin{equation}
\label{k theory of sphere}
\tilde{K}(\sphere^r) = \pi_r(BU) = \pi_{r-1}(U)
\end{equation}

A full algebraic description will be made in Section 5.

\section{Cohomology and Cohomology Ring of \texorpdfstring{$\PC^n$}{Lg}}
In this section, we will discuss the cohomology structure of $\PC^n$.
\subsection{CW Structure of \texorpdfstring{$\PC^n$}{Lg}}\mbox{}\\

$\PC^n$ has a nice CW structure:

\begin{leftbar}
\begin{proposition}
\label{cw structure of cpn}
$\PC^n$ is obtained from $\PC^{n-1}$ by attaching a single $2n$-cell.
\end{proposition}
\end{leftbar}

\begin{proof}
Let $E^{2n}$ be a 2$n$-cell.

Define $f:E^{2n} \rightarrow \PC^n$ by:
\[ f_n(z)=f_n(z_1,\cdots, z_n) = [(z_1, \cdots, z_n \sqrt{1-|z|^2})], \]
where $[(z_1, \cdots, z_n \sqrt{1-|z|^2})]$ is the equivalence class of $(z_1, \cdots, z_n \sqrt{1-|z|^2}) \in \sphere^{2n+1}$.

We first show that $f_n(\sphere^{2n-1})=\PC^{n-1}$.

Since $\sphere^{2n-1} = \lbrace z \in E^{2n} \vert \ |z|=1 \rbrace$, we have 
$f_n(\sphere^{2n-1})= \lbrace \left[(z_1, \cdots, z_n, 0) \right] \rbrace = \PC^{n-1}$.

We show that if $U^{2n}=E^{2n}-\sphere^{2n-1}$, then the restriction of $f_n$ to $U^{2n}$ 
\[\tilde{f}_n:U^{2n} \rightarrow \PC^n-\PC^{n-1} \]
is a homeomorphism.

Since $f_n(\sphere^{2n-1})=\PC^{n-1}$, we have $f_n:E^{2n} \rightarrow \PC^n$ is surjective. This implies that $f_n$ is a quotient map since $E^{2n}$ is compact and $\PC^n$ is Hausdorff.

\begin{description}
\item[$\tilde{f}_n$ is surjective] If $p \in \PC^n-\PC^{n-1}$, then $p=\left[ (z_1, \cdots, z_n, w) \right]$ with $w \not=0$ and $(z_1, \cdots, z_n, w) \in \sphere^{2n+1}$. Write $w=re^{i\theta}$, and let $\lambda= e^{-i\theta}$.
\begin{eqnarray*}
& \Rightarrow & \lambda(z_1, \cdots, z_n, w)= (\lambda z_1, \cdots, \lambda z_n, r) = (z^{'}_1, \cdots, z^{'}_n, \sqrt{1-|z|^2}) \\  & \Rightarrow & \tilde{f}_n(z^{'}) = \tilde{f} (z^{'}_1, \cdots, z^{'}_n) = (z^{'}_1, \cdots, z^{'}_n, \sqrt{1-|z|^2}) \in \left[(z_1, \cdots, z_n, w) \right]
\\ 
  & \Rightarrow & \tilde{f_n}(z^{'})=p
\end{eqnarray*}

$\therefore$ $\tilde{f}_n$ is surjective.

\item[$\tilde{f}_n$ is injective] If $\tilde{f}_n(\alpha)=\tilde{f}_n(\beta)$, then $\left[(\alpha_1, \cdots, \alpha_n, \sqrt{1-|\alpha|^2}) \right] = \left[(\beta_1, \cdots, \beta_n, \sqrt{1-|\beta|^2}) \right]$.

\begin{eqnarray*}
& \Rightarrow & (\alpha_1, \cdots, \alpha_n, \sqrt{1-|\alpha|^2})= \lambda (\beta_1, \cdots, \beta_n, \sqrt{1-|\beta|^2}), \mbox{where} \ |\lambda| = 1. 
\\
& \Rightarrow & \sqrt{1-|\alpha|^2} = \lambda \sqrt{1-|\beta|^2}
\\
& \Rightarrow & \lambda = 1
\\
& \Rightarrow & \alpha = \beta
\end{eqnarray*}
$\therefore$ $\tilde{f}_n$ is injective.
\end{description}

This concludes that $\tilde{f}_n$ is a bijective quotient map. Hence, $\tilde{f}_n$ is a homeomorphism. Moreover, since $f_n(\sphere^{2n-1})=\PC^{n-1}$, we have $\PC^n$ is obtained by $\PC^{n-1}$ by attaching a $2n$-cell. Inductively, $\PC^n$ has a CW-structure with one cell in each even dimension, and no cells of odd dimension. Therefore, the claim holds as desired.
\end{proof}

\subsection{Cohomology of \texorpdfstring{$\PC^n$}{Lg}}\mbox{}\\

We now compute the cohomology of $\PC^n$.

\begin{leftbar}
\begin{proposition}
\label{cohomology of cpn}
The cohomology of $\PC^n$ is given by:

\[ H^{2k}(\PC^n;\mathbb{Z}) = \left\{ \begin{array}{ll}
         \mathbb{Z} & \mbox{if \ $0 \leq k \leq n$};\\
        0 & \mbox{if else}.\end{array} \right. \] 

\end{proposition}
\end{leftbar}

\begin{proof}
Since  CW structure on a topological space gives a cellular filtration relative to the empty space, the $k^{th}$ cellular chain group is $\D{Z}^d$, where $d$ is the number of $k$-cells.

So combining with the CW-structure of $\PC^n$ (Proposition \ref{cw structure of cpn}), the chain complex has form:

\begin{center}
\begin{tikzpicture}[start chain] {
    \node[on chain] {$\cdots$};
    \node[on chain] {$0$};
    \node[on chain] {$0$};
    \node[on chain, label=below:{$2k$}] {$\mathbb{Z}$};
    \node[on chain] {$0$} ;
    \node[on chain] {$\mathbb{Z}$};
    \node[on chain] {$\cdots$};
    \node[on chain] {$\D{Z}$};
    \node[on chain, label=below:{$1$}] {$0$};
    \node[on chain, label=below:{$0$}] {$\mathbb{Z}$}; }
\end{tikzpicture}
\end{center}
that is, the chain complex of $\PC^n$ is:

\begin{eqnarray*}
C_{2k}(\PC^n) &=& \left\{ \begin{array}{ll}
         \mathbb{Z} & \mbox{if \ $0 \leq k \leq n$};\\
        0 & \mbox{if else}.\end{array} \right. \\
        \Rightarrow C^{2k}(\PC^n;\mathbb{Z}) &=& \mbox{Hom}(C_{2k}(\PC^n), \mathbb{Z}) \\
        &=& \left\{ \begin{array}{ll}
         \mbox{Hom}(\D{Z},\D{Z}) & \mbox{if \ $0 \leq k \leq n$};\\
        \mbox{Hom}(0, \D{Z}) & \mbox{if else}.\end{array} \right. \\
        &=& \left\{ \begin{array}{ll}
         \mathbb{Z} & \mbox{if \ $0 \leq k \leq n$};\\
        0 & \mbox{if else}.\end{array} \right.
        \end{eqnarray*}
        
$\therefore$ $H^{2k}(\PC^n;\mathbb{Z}) = \left\{ \begin{array}{ll}
         \mathbb{Z} & \mbox{if \ $0 \leq k \leq n$};\\
        0 & \mbox{if else}.\end{array} \right. $

\end{proof}

\subsection{Cohomology Ring of \texorpdfstring{$\PC^n$}{Lg}}\mbox{}\\

To compute the cohomology ring of $\PC^n$, we need the following lemma:

\begin{leftbar}
\begin{lemma}[{{\cite[{\bf p.152 Corollary}]{may}}}]
\label{lemma for cohomology ring}
Let $T_p \subset H^p(M)$ be the torsion subgroup. The cup product pairing

\[ \alpha \otimes \beta \rightarrow\langle \alpha \beta ,z \rangle \]

induces a non-singular pairing

\[ H^p(M;\D{Z})/T_p \otimes H^{n-p}(M;\D{Z})/T_{n-p} \rightarrow \mathbb{Z} \]
\end{lemma}
\end{leftbar}

\begin{proof}
Recall that if $M$ is a compact $n$-manifold, then $H_q(M;\D{Z})$ is finitely generated for all $q$, and $H^n(M;\D{Z}) = \D{Z}$.

If $\alpha \in T_p$, then there exists $r \in \D{Z}$ such that $r\alpha = 0$. For $\beta \in H^{n-p}(M;\D{Z})$, we have $r(\alpha \cup \beta) = 0$. So $\alpha \cup \beta = 0$ since $H^n(M;\D{Z}) = \D{Z}$. This shows paring vanishes on torsion elements.

Now, since $\Ext^1_{\D{Z}}(\D{Z}_r, \D{Z}) = \D{Z}_r$, and each $H_p(M;\D{Z})$ is finitely generated, it follows that $\Ext^1_{\D{Z}}(M;\D{Z})$ is a torsion group. By Universal Coefficient Theorem, we have:

\[ H^p(M;\D{Z})/T_p = \mbox{Hom}(H_p(M;\D{Z}),\D{Z}).\] 

Hence, if $\alpha \in H^p(M;\D{Z})$ projects to a generator of the free abelian group $H_p(M;\D{Z})/T_p$, then there exists $a \in H_p(M;\D{Z})$ such that $\langle \alpha, a \rangle = 1$.

By Poincar\'e Duality, there exists $\beta \in H^{n-p}(M;\D{Z})$ such that $\beta \cap z = a$, where $z$ is a generator of $H^n(M;\D{Z})$. So $\langle \beta \cup \alpha, z \rangle = \langle \alpha, \beta \cap z \rangle = 1$. The claim holds as desired.
\end{proof}

\begin{leftbar}
\begin{proposition}
\label{cohomology ring of cpn}
The cohomology ring of $\PC^n$ is given by:

\[ H^{\ast}(\PC^n;\mathbb{Z}) = \mathbb{Z}[\alpha]/\langle \alpha^{n+1}\rangle \]

where $\alpha$ is a generator of $H^2(\PC^n; \D{Z})$.
\end{proposition}
\end{leftbar}

\begin{proof}
We will prove the claim by induction on $n$.

When $n=1$, $\PC^1$ is homeomorphic to $\sphere^2$ by the CW-structure of $\PC^n$ (Proposition \ref{cw structure of cpn}). Since $H^1(\sphere^2)=0$, we have:
\begin{eqnarray*}
H^\ast(\PC^1;\D{Z}) & = & H^0(\PC^1;\D{Z}) \oplus H^2(\PC^1;\D{Z}) \\
             & = & \D{Z} \oplus \D{Z} \\
             &=& \D{Z}[\alpha]/\langle \alpha^2 \rangle \ \mbox{since $\alpha^2 \in H^4(\PC^1;\D{Z}) = 0$}            
\end{eqnarray*}

So the claim holds for $n=1$.

Assume the claim holds for $n\in \lbrace 1, \cdots, k-1 \rbrace \subseteq \D{N}$. When $n=k$, the assumption asserts that if $\alpha$ generates $H^2(\PC^n;\D{Z})$, then $\alpha^q$ generates $H^{2q}(\PC^n;\D{Z})$ for $q<n$. By Lemma \ref{lemma for cohomology ring}, there exists $\beta \in H^{2n-2}(\PC^n;\D{Z})$ such that 
\[\langle \alpha \cup \beta, z \rangle = 1 \]
where $z$ is a generator of $H^{2n}(\PC^n;\D{Z})$
Note that $\beta$ must be a generator, so we have $\beta = \pm \alpha^{n-1}$. Hence, $\alpha^n$ generates $H^{2n}(\PC^n;\D{Z})$, and $\alpha^{n+1}=0$ since $\alpha^{n+1} \in H^{2n+2}(\PC^n;\D{Z})=0$.

$\therefore$ By induction, the claim holds as desired.
\end{proof}

\section{Chern Character}
In this section, we will introduce the basic notations and concepts about Chern character, and compute the Chern character of the Hopf bundle. For further details of the Chern theory, see \cite{milnorsta} or \cite{madsen}.

\subsection{Gysin Sequence and Chern Classes}\mbox{}\\

Suppose $\omega$ is complex $n$-dimensional vector bundle, with total space $E$ and base space $B$. We construct a canonical $(n-1)$-dimensional vector bundle $\omega_0$ over the deleted total space $E_0$ (we obtain $E_0$ by deleting the zero section of $E$). A point in $E_0$ is determined by a fiber $F$ of $\omega$, and a non-zero vector $v \in F$. We define the fiber of $\omega_0$ over $v$ to be the quotient vector space $F/(\D{C}v)$. By $\D{C}v$ we mean the 1-dimensional subspace spanned by the non-zero vector $v$. From this construction, $F/(\D{C}v)$ is a complex vector space of dimension $n-1$, and clearly can be considered as fiber of $\omega_0$. So $\omega_0$ is indeed a vector bundle over $E_0$.

By \cite[p.~143]{milnorsta}, any real oriented $2n$-plane bundle $\pi:E \rightarrow B$ possesses an exact {\bf Gysin sequence} with integer coefficients:
\begin{center}
\begin{tikzpicture}[start chain] {
    \node[on chain] {$\cdots$};
    \node[on chain] {$H^{i-2n}(B)$};
    \node[on chain, join={node[above] {$\cup_c$}}] {$H^i(B)$};
    \node[on chain, join={node[above] {$\pi^{\ast}_0$}}] {$H^i(E_0)$};
    \node[on chain] {$H^{i-2n+1}(B)$} ;
    \node[on chain] {$\cdots$};
        }
\end{tikzpicture}
\end{center}
where 
\[ \cup_c: H^i(E) \rightarrow H^{i+n}(E,E_0) \ \mbox{with} \ \cup_c(u)=u\cup c \]
is an isomorphism for any arbitrary coefficient module. 

Note that for $i<2n-1$, the groups $H^{i-2n}(B)$ and $H^{i-2n+1}(B)$ are zero, so the map:
\[ \pi^{\ast}_0: H^i(B) \rightarrow H^i(E_0)\]
 is an isomorphism. We now define the Chern class of a vector bundle.

\begin{leftbar}
\begin{definition}[{{\bf Chern classes}}]
\label{ith chern class}
Let $\omega$ be a complex $n$-plane bundle, with total space $E$ and base space $B$. Denote $\omega_0$ be the $(n-1)$-bundle over the deleted total space $E_0$. 

For $i<n$, the $\bm{i^{th} \mbox{ \bf Chern Class of} \ \omega}$ is defined to be:

\[ c_i(\omega) = \left\{ \begin{array}{ll}
         \pi^{\ast^{-1}}_0(\omega_0) & \mbox{if $0 < i < n$};\\
         1 & \mbox{if $i=0$}.\end{array} \right. \]
where $\pi^{\ast}_0$ is defined in the Gysin sequence:

\begin{center}
\begin{tikzpicture}[start chain] {
    \node[on chain] {$\cdots$};
    \node[on chain] {$H^{2i-2n}(B)$};
    \node[on chain, join={node[above] {$\cup_c$}}] {$H^{2i}(B)$};
    \node[on chain, join={node[above] {$\pi^{\ast}_0$}}] {$H^{2i}(E_0)$};
    \node[on chain] {$H^{2i-2n+1}(B)$} ;
    \node[on chain] {$\cdots$};
        }
\end{tikzpicture}
\end{center}
Note that $c_i(\omega) \in H^{2i}(B;\D{Z})$.

For $i > n$ we just set $c_i(\omega) = 0$. 

If $i=n$, the {\bf top Chern class} is defined to be the Euler class:
Put $i=n$, the Gysin sequence gives:
\begin{center}
\begin{tikzpicture}[start chain] {
    \node[on chain] {$H^0(B)$};
    \node[on chain] {$H^{2n}(B)$};
    \node[on chain] {$H^{2n}(E_0)$};
        }
\end{tikzpicture}
\end{center}
which the composition of maps sends $1 \in H^0(B)$ to some element $y \in H^{2n}(E_0)$. Then, $c_n(\omega)$ is the (unique) element in $H^{2n}(B)$ that is the pre-image of $y$.

The {\bf total Chern class} of $E$ is the formal sum of Chern classes:
\[ c(E)=1+c_1(E)+c_2(E)+\cdots+c_n(E) \]
which is an element in the even cohomology $H^{\mbox{even}}(B)= \displaystyle\bigoplus_{i=0}^nH^{2i}(B)$.
\end{definition}
\end{leftbar}

\begin{description}

\item[Note] $\pi^{\ast}_0:H^{2i}(B) \rightarrow H^{2i}(E_0)$ is an isomorphism for $i<n$. So the $i^{th}$ Chern class is well-defined and unique.
\end{description}

\begin{leftbar}
\begin{definition}[{{\bf Chern Character}}]
\label{chern character}
Let $x_1, \cdots, x_n$ be variables, and $p_k(x_1, \cdots, x_n)$ be the $k^{th}$ power sum:
\[ p_k(x_1, \cdots, x_n) = \displaystyle \sum_{i=1}^n x_i^k = x_1^k + \cdots + x_n^k. \]

Moreover, denote $e_i$ be the elementary symmetric polynomial:
\begin{eqnarray*}
 e_0(x_1, \cdots, x_n) &=& 1; \\
 e_1(x_1, \cdots, x_n) &=& x_1 + \cdots + x_n; \\
 e_2(x_1, \cdots, x_n) &=& \displaystyle \sum_{1 \leq i \leq j \leq n} x_i x_j; \\
 e_n(x_1, \cdots, x_n) &=& x_1 x_2 \cdots x_n; \\
 e_k(x_1, \cdots, x_n) &=& 0 \ \mbox{for $k>n$}.
\end{eqnarray*}

The functions $p_k$ can be expressed as polynomials of $e_1, \cdots, e_n$. For example:
\begin{eqnarray*}
p_1 &=& e_1; \\
p_2 &=& e_1^2 - 2e_2;\\
p_3 &=& e_1^3 -3e_1 e_2 + 3e_3; \\
\vdots
\end{eqnarray*}
As a result, we can write $p_k = s_k (e_1, \cdots, e_n)$.

Let $E$ be a vector bundle of rank $n$. The  $\bm{k^{th}}$ {\bf Chern character} of $E$ is:
\[ ch_k(E) = \frac{s_k(c(E))}{k!} \]
for $k > 0$. 

If $k=0$, we define:
\[ ch_0(E) = \mbox{dim}(E). \]

The {\bf total Chern character} of $E$ is defined to be:

\begin{equation}
\label{ch(E)}
ch(E) = \mbox{dim}(E) + \displaystyle\sum_{k=1}^{\infty} \frac{s_k(c(E))}{k!}
      = n + c_1(E) + \frac{c_1^2(E)-2c_2(E)}{2} +  \frac{c_1^3(E)-3c_1(E)c_2(E)+3c_3(E)}{3!} + \cdots
\end{equation}
\end{definition}
\end{leftbar}

We list some properties of the Chern character and Chern class. For the proofs, see \cite{milnorsta} or \cite{madsen}.

\begin{leftbar}
\begin{proposition}
\label{sum and product of ch}
Suppose $E_0$, $E_1$ are two complex vector bundles over a common paracompact base space $B$, then:
\begin{eqnarray}
ch(E_0 \oplus E_1) &=& ch(E_0) + ch(E_1);\\
ch(E_0 \otimes E_1) &=& ch(E_0)ch(E_1)
\end{eqnarray}
(The product on the right hand side is referring to the cup product).
Moreover, we have:
\[ c(E_0 \oplus E_1) = c(E_0)c(E_1). \]
\end{proposition}
\end{leftbar}

The additivity of Chern character allows us to extend the Chern character to a homomorphism in $K$-theory:
\[ \mbox{ch}:K(X) \rightarrow \mbox{H}^{\mbox{even}}(X;\D{Q}).\]
Equally important, the multiplicative proper of Chern character makes this homomorphism multiplicative (and hence becomes a ring homomorphism). As a result, the Chern character contributes on investigating the algebraic properties of elements of the $K$-group.

The following proposition gives an application of the Chern character to $K$-theory:
\begin{leftbar} 
\begin{proposition}
\label{linearly independent}
Suppose $E_0$, $E_1$ are two complex vector bundles over a common paracompact base space $B$. If $ch(E_0) \not= ch(E_1)$, then $E_0$ and $E_1$ give different elements in the $K$-group.
\end{proposition}
\end{leftbar}

\subsection{Chern Character of Hopf Bundle}\mbox{}\\

Recall that the Hopf bundle $\zeta$ is a rank $1$ bundle $\sphere^{2n+1} \rightarrow \PC^n$ over $\PC^n$. Here, we compute the Chern character of $\zeta$. In this section, the Chern character of the Hopf bundle is computed.

\begin{leftbar}
\begin{lemma}
\label{c_i zeta}
For $2 \leq i \leq n$, $c_i(\zeta)=0$.
\end{lemma}
\end{leftbar}

\begin{proof}
Suppose $n=2$. Since $\zeta$ has rank $1$, and $n=2>1$, so by Definition~\ref{ith chern class}, we immediately have $c_2(\zeta) =0$. Therefore, the claim holds for $n=2$.

Assume the claim holds for $n \in \lbrace 2, \cdots, k-1 \rbrace \subseteq \D{N}$.
If $n=k$, then by definition:
\[ c_i(\zeta) = \pi_0^{\ast^{-1}}(\zeta_0),\]
where $\zeta_0$ is a bundle over $\sphere^{2n}$. By assumption, $c_i(\zeta_0) \in H^{2i}(\sphere^{2n-1}) = 0$ for every $2 \leq i \leq n$. Since $\pi_0^{\ast}$ is an isomorphism, we must have $c_i(\zeta)= \pi_0^{\ast^{-1}}(c_i(\zeta_0)) = \pi_0^{\ast^{-1}}(0) = 0$.

$\therefore$ By induction, the claim holds as desired.
\end{proof}

\begin{leftbar}
\begin{proposition}
Let $x$ be a generator of $H^2(\PC^n; \D{Z})$. Then $ch(\zeta)=e^x$.
\end{proposition}
\end{leftbar}

\begin{proof}
Using Equation \ref{ch(E)}, we have:
\begin{eqnarray*}
ch(\zeta) &=& 1+ c_1(\zeta) + \frac{c_1^2(\zeta)-2c_2(\zeta)}{2!} + \frac{c_1^3(\zeta)-3c_1(\zeta)c_2(\zeta)+3c_3(\zeta)}{3!} + \cdots \\
&=& 1+c_1(\zeta) + \frac{c_1^2(\zeta)}{2!} + \frac{c_1^3(\zeta)}{3!} + \cdots (\mbox{by Lemma \ref{c_i zeta}}) \\
&=& e^x
\end{eqnarray*}
\end{proof}

\section{Bott Periodicity Theorem}
It is often hard to determine the equivalence classes of a given complex vector bundles, even we have a natural identification in Theorem~\ref{k theory and bu}:
\[ \tilde{K}(X) = \pi (X, BU). \]
However, the {\bf Bott Periodicity Theorem} provides a powerful machinery to compute $K$-groups. If $X=\D{C}$, then we only need to compute $K^0$ and $K^1$. If $X=\D{R}$, we need to compute $K^0, K^1, \cdots, K^7$. As a result, we are able to extend the group $K^q(X,A)$ to all integers $q$. In this section, we will give the statement of Bott Periodicity Theorem and two of its corollaries. For the proofs, see \cite{fuchs}.

\begin{leftbar}
\begin{theorem}[{{\bf Bott Periodicity}}]
\label{bott}
For any finite CW complex, the Bott map: 
\[ \mbox{Bott}: K(X) \oplus K(X) \rightarrow K(X \times \sphere^2)\] 
defined by
\[ \mbox{Bott}(\alpha_1, \alpha_2) = (\alpha_1 \otimes \tilde{1}) + (\alpha_2 \otimes \zeta)\]
is an isomorphism.

Here, $\tilde{1}$ is the trivial bundle of dimension $1$, and $\zeta$ is the Hopf bundle.
\end{theorem}
\end{leftbar}

An immediately consequence of Bott Periodicity Theorem is that:
\begin{leftbar}
\begin{corollary}
\label{1st corollary of Bott}
Let $X$ be a finite CW complex. Then
\[ \tilde{K}(X) \cong \tilde{K}(\Sigma^2X) \]
\end{corollary}
\end{leftbar}

Because of Bott Periodicity Theorem, to compute the $K$-theory of complex vector bundle, we only need to calculate $K^0$ and $K^1$. Moreover, by elementary homotopy theory we have:
\[ \pi_1(U) = \D{Z} \ \mbox{and} \ \pi_1(BU) = 0. \]
So we have the following corollary, which calculates the $K$-theory of sphere.

\begin{leftbar}
\begin{corollary}[{{\bf Corollary of Corollary \ref{1st corollary of Bott}}}]
\label{reduced k theory of sphere}
For any $q$,
\[ \pi_{i-2}(BU) \cong \pi_i(BU) = \tilde{K}(\sphere^i) =  \left\{ \begin{array}{ll}
         \D{Z} & \mbox{if $i$ is even};\\
        0 & \mbox{if $i$ is odd}.\end{array} \right. 
\]
and hence,
\[ \pi_{i-2}(U) \cong \pi_i(U) = \left\{ \begin{array}{ll}
         0 & \mbox{if $i$ is even};\\
        \D{Z} & \mbox{if $i$ is odd}.\end{array} \right.  \]
\end{corollary}
\end{leftbar}

So far, the group $K^q(X)$ is defined for $q<0$. Because of Corollary~\ref{reduced k theory of sphere}, we can write:
\[ K^q(X) \cong K^{q+2}(X) \]
for all $q$. Therefore, the group $K^q(X)$ is defined for every integer $q$. This completes the construction of $K$-groups.

\section{Proof of Main Theorem}
To prove Theorem \ref{hahaj}, we still need one lemma:

\begin{leftbar}
\begin{lemma}
\label{image of ch}
The homomorphism $\mbox{ch}: \tilde{K}(\sphere^{2n}) \rightarrow \tilde{H}^{\ast}(\sphere^{2n};\D{Q})$ has image $\D{Z}$.
\end{lemma}
\end{leftbar}

\begin{proof}
When $n=1$, we have $\tilde{K}(\sphere^2) = \D{Z}$ and $\tilde{H}^{\ast}(\sphere^2;\D{Q}) = H^{2n}(\sphere^2; \D{Q}) = \D{Q}$. The claim follows immediately, since ch maps the trivial bundle $1$ to the identity element in $\D{Q}$ (which is the natural number $1$).

For general $n$, just observe that the following diagram commutes (by \cite{fuchs}):
\begin{center}
\begin{tikzpicture}
  \node (A) {$\tilde{K}^0(X)$}; 
  \node (B) [below=of A] {$\tilde{K}^2(X)$}; 
  \node (C) [right=of A] {$\tilde{H}^{\mbox{even}}(X;\D{Q})$}; 
  \node (D) [right=of B] {$\tilde{H}^{\mbox{even}}(X;\D{Q})$};
  \node (E) [below=of B] {$\tilde{K}^0(X)$};
  \node (F) [right=of E] {$\tilde{H}^{\mbox{even}}(X;\D{Q})$};

  \draw[-stealth] (A)-- node[left] {\small $\Sigma^2$} (B); 
  \draw[-stealth] (B)-- node [above] {\small ch} (D); 
  \draw[-stealth] (A)-- node [above] {\small ch} (C); 
  \draw[-stealth] (C)-- node [right] {\small $\Sigma^2$} (D); 
  \draw[-stealth] (E)-- node [above] {\small ch} (F);
  \draw[-stealth] (B)-- node [left] {\small Bott} (E);
  \draw[white] (D) edge node [rotate=90, black] {$=$} (F);
\end{tikzpicture}
\end{center}
So the claim holds.
\end{proof}
Let's restate the main theorem first:

\noindent\usebox{\maintheorem}

\begin{proof}
By Bott Periodicity Theorem (Theorem \ref{bott}), we can focus on the cases $q=0$ and $q=1$. 

By the CW structure of $\PC^n$ (Proposition \ref{cw structure of cpn}), we have $\PC^n / \PC^{n-1}$ is homeomorphic to $\sphere^{2n}$. We will compute $K^0(\PC^n)$ by induction on $n$.

If $n=1$, $\PC^n$ is homeomorphic to $\sphere^2$.
\begin{eqnarray*}
\Rightarrow K^0(\PC^n) & = & K^0(\sphere^2) \ (\mbox{since $K(X)$ is homotopy invariant.}) \\
             & = & \tilde{K}(\sphere^2) \oplus \D{Z} \ (\mbox{by Definition \ref{def of higher K-groups and relative K-groups}})\\
             & = & \D{Z} \oplus \D{Z} \ (\mbox{by Equation \ref{k theory of sphere}})
\end{eqnarray*}
and
\begin{eqnarray*}
K^1(\PC^n) & = & K^1(\sphere^2)\\
             & = & \tilde{K}(\sphere^3) \ (\mbox{by Definition \ref{def of higher K-groups and relative K-groups}})\\
             & = & \pi_3(BU) \ (\mbox{by Corollary \ref{reduced k theory of sphere}}) \\
             &=& 0
\end{eqnarray*}

So the claim holds for $n=1$. Assume the claim holds for $n \in \lbrace 1, \cdots, k \rbrace \subseteq \D{N}$. If $n=k+1$, by the CW pair $(\PC^{k+1}, \PC^{k})$, the part of exact sequence of $\tilde{K}$-groups yields:
\begin{center}
\begin{tikzpicture}
  \node (A) {$0$}; 
  \node (B) [right= of A] {$\tilde{K}(\PC^k)$}; 
  \node (C) [right= of B] {$\tilde{K}(\PC^{k+1})$}; 
  \node (D) [right= of C] {$\tilde{K}(\PC^{k+1}/ \PC^k)$};
  \node (E) [right= of D] {$\tilde{K}(\sum \PC^k)$};

  \draw[stealth-] (A)-- node[right] {\small} (B); 
  \draw[stealth-] (B)-- node [right] {\small} (C); 
  \draw[stealth-] (C)-- node [right] {\small} (D);
  \draw[stealth-] (D)-- node [right] {\small} (E);
    
\end{tikzpicture}
\end{center}
which reduced to the following diagram:
\begin{center}
\begin{tikzpicture}
  \node (A) at(-5,0) {$0$}; 
  \node (B) at(-3,0) {$\tilde{K}(\PC^k)$}; 
  \node (C) at(0,0) {$\tilde{K}(\PC^{k+1})$}; 
  \node (D) at(3,0) {$\tilde{K}(\PC^{k+1}/ \PC^k)$};
  \node (E) at(6,0) {$\tilde{K}(\sum \PC^k)$};
  \node (F) at(-5,-1) {$0$};
  \node (G) at(-3,-1) {$\D{Z}^{\oplus k}$};
  \node (H) at(0,-1) {$\tilde{K}(\PC^{k+1})$};
  \node (I) at(3,-1) {$\tilde{K}(\sphere^{2k+2})$};
  \node (J) at(6,-1) {$\tilde{K}(\sum \PC^k)$};
  \node (K) at(-5,-2) {$0$};
  \node (L) at(-3,-2) {$\D{Z}^{\oplus k}$};
  \node (M) at(0,-2) {$\tilde{K}(\PC^{k+1})$};
  \node (N) at(3,-2) {$\D{Z}$};
  \node (O) at(6,-2) {$\tilde{K}(\sum \PC^k)$};

  \draw[stealth-] (A)-- node [right] {\small} (B); 
  \draw[stealth-] (B)-- node [right] {\small} (C); 
  \draw[stealth-] (C)-- node [right] {\small} (D);
  \draw[stealth-] (D)-- node [right] {\small} (E);
  \draw[stealth-] (F)-- node [right] {\small} (G);
  \draw[stealth-] (G)-- node [right] {\small} (H);
  \draw[stealth-] (H)-- node [right] {\small} (I);
  \draw[stealth-] (I)-- node [right] {\small} (J);
  \draw[stealth-] (K)-- node [right] {\small} (L);
  \draw[stealth-] (L)-- node [right] {\small} (M);
  \draw[stealth-] (M)-- node [right] {\small} (N);
  \draw[stealth-] (N)-- node [right] {\small} (O);
  \draw[white] (A)  edge node [rotate=90, black] {$=$} (F);
  \draw[white] (B)  edge node [rotate=90, black] {$=$} (G);
  \draw[white] (C)  edge node [rotate=90, black] {$=$} (H);
  \draw[white] (D)  edge node [rotate=90, black] {$=$} (I);
  \draw[white] (E)  edge node [rotate=90, black] {$=$} (J);
  \draw[white] (F)  edge node [rotate=90, black] {$=$} (K);
  \draw[white] (G)  edge node [rotate=90, black] {$=$} (L);
  \draw[white] (H)  edge node [rotate=90, black] {$=$} (M);
  \draw[white] (I)  edge node [rotate=90, black] {$=$} (N);
  \draw[white] (J)  edge node [rotate=90, black] {$=$} (O);

\end{tikzpicture}
\end{center}

By assumption, $\tilde{K}(\PC^k) = \D{Z}^{\oplus k}$. By Proposition \ref{cw structure of cpn}, we have $\sphere^{2k}=\PC^{k+1}/\PC^k$. So we obtain the second row from the first row. By Corollary \ref{reduced k theory of sphere}, we have $\tilde{K}(\sphere^{2k+2}) = \D{Z}$, so we obtain the third row.

Now, apply Five Lemma to the diagram:
\begin{center}
\begin{tikzpicture}
  \node (K) at(-4,0) {$0$};
  \node (L) at(-2,0) {$\D{Z}^{\oplus k}$};
  \node (M) at(0,0) {$\tilde{K}(\PC^{k+1})$};
  \node (N) at(2,0) {$\D{Z}$};
  \node (O) at(4,0) {$\tilde{K}(\sum \PC^k)$};
  
  \node (A) at(-4,-1) {$0$};
  \node (B) at(-2,-1) {$\D{Z}^{\oplus k}$};
  \node (C) at(0,-1) {$\D{Z} \oplus \D{Z}^{\oplus k}$};
  \node (D) at(2,-1) {$\D{Z}$};
  \node (E) at(4,-1) {$\tilde{K}(\sum \PC^k)$};
  
  \draw[stealth-] (K)-- node [right] {\small} (L);
  \draw[stealth-] (L)-- node [right] {\small} (M);
  \draw[stealth-] (M)-- node [right] {\small} (N);
  \draw[stealth-] (N)-- node [right] {\small} (O);
  \draw[stealth-] (A)-- node [right] {\small} (B);
  \draw[stealth-] (B)-- node [right] {\small} (C);
  \draw[stealth-] (C)-- node [right] {\small} (D);
  \draw[stealth-] (D)-- node [right] {\small} (E);
  
  \draw[white] (K)  edge node [rotate=90, black] {$=$} (A);
  \draw[white] (L)  edge node [rotate=90, black] {$=$} (B);
  \draw[white] (N)  edge node [rotate=90, black] {$=$} (D);
  \draw[white] (O)  edge node [rotate=90, black] {$=$} (E);
  \draw[-stealth] (M)-- node [right] {\small} (C);
\end{tikzpicture}
\end{center}

We have $\tilde{K}(\PC^{k+1}) = \D{Z} \oplus \D{Z}^{\oplus k} = \D{Z}^{\oplus k+1}$. Since $K(X) = \tilde{K}(X) \oplus \D{Z}$, we have \[K^0(\PC^{k+1}) = \D{Z}^{\oplus k+1} \oplus \D{Z} = \D{Z}^{\oplus k+2}.\]

We do the same thing to get $K^1(\PC^n)$. By assuming $K^1(\PC^n) = 0$ for $n \in \lbrace 1, \cdots, k \rbrace \subseteq \D{N}$, when $n=k+1$, we have the exact sequence:

\begin{center}
\begin{tikzpicture}
  \node (A) at(-6,0) {$K^0(\PC^{k+1},\PC^k)$}; 
  \node (B) at(-3,0) {$K^{-1}(\PC^k)$}; 
  \node (C) at(0,0) {$K^{-1}(\PC^{k+1})$}; 
  \node (D) at(4,0) {$K^{-1}(\PC^{k+1}, \PC^k)$};
  \node (E) at(8,0) {$K^{-2}(\PC^k)$};
  \node (F) at(-6,-1) {$K^0(\PC^{k+1},\PC^k)$};
  \node (G) at(-3,-1) {$K^1(\PC^k)$};
  \node (H) at(0,-1) {$K^1(\PC^{k+1})$};
  \node (I) at(4,-1) {$K^1(\PC^{k+1}, \PC^k)$};
  \node (J) at(8,-1) {$K(\PC^k)$};
  \node (K) at(-6,-2) {$\D{Z}$};
  \node (L) at(-3,-2) {$0$};
  \node (M) at(0,-2) {$K^1(\PC^k)$};
  \node (N) at(4,-2) {$0$};
  \node (O) at(8,-2) {$\D{Z}^{\oplus k+1}$};

  \draw[stealth-] (A)-- node [right] {\small} (B); 
  \draw[stealth-] (B)-- node [right] {\small} (C); 
  \draw[stealth-] (C)-- node [right] {\small} (D);
  \draw[stealth-] (D)-- node [right] {\small} (E);
  \draw[stealth-] (F)-- node [right] {\small} (G);
  \draw[stealth-] (G)-- node [right] {\small} (H);
  \draw[stealth-] (H)-- node [right] {\small} (I);
  \draw[stealth-] (I)-- node [right] {\small} (J);
  \draw[stealth-] (K)-- node [right] {\small} (L);
  \draw[stealth-] (L)-- node [right] {\small} (M);
  \draw[stealth-] (M)-- node [right] {\small} (N);
  \draw[stealth-] (N)-- node [right] {\small} (O);
  \draw[white] (A)  edge node [rotate=90, black] {$=$} (F);
  \draw[white] (B)  edge node [rotate=90, black] {$=$} (G);
  \draw[white] (C)  edge node [rotate=90, black] {$=$} (H);
  \draw[white] (D)  edge node [rotate=90, black] {$=$} (I);
  \draw[white] (E)  edge node [rotate=90, black] {$=$} (J);
  \draw[white] (F)  edge node [rotate=90, black] {$=$} (K);
  \draw[white] (G)  edge node [rotate=90, black] {$=$} (L);
  \draw[white] (H)  edge node [rotate=90, black] {$=$} (M);
  \draw[white] (I)  edge node [rotate=90, black] {$=$} (N);
  \draw[white] (J)  edge node [rotate=90, black] {$=$} (O);

\end{tikzpicture}
\end{center}
By Bott Periodicity Theorem, we have $K^{-1}(X) = K^1(X)$, so we get the second row from the first row. By assumption, we have $K^1(\PC^k) = 0$. Moreover, we have:
\[ K^1(\PC^{k+1}, \PC^k) = K^1(\sphere^{2k+2}) =  \tilde{K}(\sphere^{2k+3}) = 0 \ \mbox{(by Corollary \ref{reduced k theory of sphere}).} \]
Combining the result we get for $K^0(\PC^n)$, we get the third row from the second row.

Again, apply the Five Lemma to the diagram:
\begin{center}
\begin{tikzpicture}
  \node (K) at(-4,0) {$\D{Z}$};
  \node (L) at(-2,0) {$0$};
  \node (M) at(0,0) {$K^1(\PC^{k+1})$};
  \node (N) at(2,0) {$0$};
  \node (O) at(4,0) {$\D{Z}^{\oplus k+1}$};
  
  \node (A) at(-4,-1) {$\D{Z}$};
  \node (B) at(-2,-1) {$0$};
  \node (C) at(0,-1) {$0$};
  \node (D) at(2,-1) {$0$};
  \node (E) at(4,-1) {$\D{Z}^{\oplus k+1}$};
  
  \draw[stealth-] (K)-- node [right] {\small} (L);
  \draw[stealth-] (L)-- node [right] {\small} (M);
  \draw[stealth-] (M)-- node [right] {\small} (N);
  \draw[stealth-] (N)-- node [right] {\small} (O);
  \draw[stealth-] (A)-- node [right] {\small} (B);
  \draw[stealth-] (B)-- node [right] {\small} (C);
  \draw[stealth-] (C)-- node [right] {\small} (D);
  \draw[stealth-] (D)-- node [right] {\small} (E);
  
  \draw[white] (K)  edge node [rotate=90, black] {$=$} (A);
  \draw[white] (L)  edge node [rotate=90, black] {$=$} (B);
  \draw[white] (N)  edge node [rotate=90, black] {$=$} (D);
  \draw[white] (O)  edge node [rotate=90, black] {$=$} (E);
  \draw[-stealth] (M)-- node [right] {\small} (C);
\end{tikzpicture}
\end{center}
We have $K^1(\PC^{k+1}) = 0$. So by induction, the claim holds as desired.

For the ring structure of $K(\PC^n)$, we first compute the Chern character of $\gamma$.
By Proposition \ref{c_i zeta}, we have:
\[ \mbox{ch}(\zeta) = e^x = 1 + x + \frac{x^2}{2} + \cdots \]
where $x$ is a generator of $H^2(\PC^n; \D{Z})$.

Using Proposition \ref{sum and product of ch}, we get:
\[ \mbox{ch}(\gamma) = \mbox{ch}(\zeta - \tilde{1}) = e^x - 1 = x + \frac{x^2}{2} + \cdots, \]
and
\[ \mbox{ch}(\gamma^k) = \mbox{ch}(\gamma^{\otimes k}) = (\mbox{ch}(\gamma))^k = x^k + \frac{k}{2}x^{k+1} + \cdots. \]
for $1<k<n$.

Now, when $k=n+1$, by Proposition \ref{cohomology ring of cpn}, we have $x^{n+1} = 0$. So we get
\[ \mbox{ch}(\gamma^{n+1}) = x^{n+1} = 0 \]
and 
\[ \mbox{ch}(\gamma^n) = x^n. \]

By Proposition \ref{linearly independent}, $\tilde{1}, \gamma, \cdots, \gamma^n$ represent different elements in $K(\PC^n)$. We show that they generate $K^{\ast}(\PC^n)$ over $\D{Z}$ by induction on $n$.

When $n=1$, the claim follows immediately. So assume the claim holds for $n \in \lbrace 1, \cdots, k-1 \rbrace \subseteq \D{N}$.

If $n = k$, the part of exact sequence of $K$-groups yields:

\begin{center}
\begin{tikzpicture}
  \node (A) {$K(\PC^n, \PC^{n-1}) = \tilde{K}(\sphere^{2n})$}; 
  \node (B) [right= of A] {$K(\PC^n)$}; 
  \node (C) [right= of B] {$K(\PC^{n-1})$}; 
  \node (D) [right= of C] {$0.$};

  \draw[-stealth] (A)-- node[above] {\small $p^{\ast}$} (B); 
  \draw[-stealth] (B)-- node [above] {\small $i^{\ast}$} (C); 
  \draw[-stealth] (C)-- node [right] {\small} (D);
   
\end{tikzpicture}
\end{center}

By assumption, $K(\PC^{n-1})$ is generated by $\tilde{1}, \gamma, \cdots, \gamma^{n-1}$. Consider a map:
\[ g: K(\PC^{n-1}) \rightarrow K(\PC^n) \]
such that $g(\gamma) = \gamma$. It is clear that $i^{\ast} \circ g$ is the identity map on $K(\PC^{n-1})$. So the previous exact sequence of $K$-groups splits:

\[ K(\PC^n) = \mbox{Im}(p^{\ast}) \oplus \mbox{Im}(q). \]

Thus, if $\alpha \in K(\PC^n)$, then $\alpha$ can be written as:

\begin{eqnarray}
\label{alpha}
\nonumber \alpha &=& g \circ i^{\ast} (\alpha) + p^{\ast}(\beta)\\
&=& r_0 + r_1 \gamma + \cdots + r_{n-1} \gamma^{n-1} + p^{\ast}(\beta),
\end{eqnarray}
where $r_i \in \D{Z}$, and $\beta \in \tilde{K}(\sphere^{2n})$.

We wish to show that $p^{\ast}(\beta) = r \gamma^n$ for some $r \in \D{Z}$, which can be done by computing the Chern character of $p^{\ast}(\beta))$. We have the following commutative diagram:

\begin{center}
\begin{tikzpicture}
  \node (A) at(0,0) {$\tilde{K}(\sphere^{2n})$}; 
  \node (B) at(0,-2) {$H^{2n}(\sphere^{2n})$}; 
  \node (C) at(4,0) {$K(\PC^n)$}; 
  \node (D) at(4,-2) {$H^{2n}(\PC^n)$};
   
  \draw[-stealth] (A)-- node[left] {\small ch} (B); 
  \draw[-stealth] (B)-- node [above] {\small $p^{\ast}$} (D); 
  \draw[-stealth] (A)-- node [above] {\small $p^{\ast}$} (C); 
  \draw[-stealth] (C)-- node [right] {\small ch} (D); 
\end{tikzpicture}
\end{center}

Furthermore, by Proposition~\ref{cohomology of cpn}, we know that
\[ H^{2n}(\PC^n;\D{Z}) = \D{Z} = \D{Z}x^n \]
where $x^n = \mbox{ch}(\gamma^n)$.

So we must have $\mbox{ch}(p^{\ast}(\beta)) = r x^n$ for some $r \in \D{Z}$.

Now, since $H^{2n}(\PC^n;\D{Z}) = \D{Z}$, it is torsion free. So the map:
\[  \mbox{ch}: K(\PC^n) \rightarrow H^{2n}(\PC^n) \]
has trivial kernel. Combining this with the fact that $\mbox{ch}(\gamma^n) = x^n$, we must have:
\[ p^{\ast}(\beta) = r \gamma^n. \]

Therefore, Equation~\ref{alpha} becomes:

\begin{eqnarray*}
\alpha &=& g \circ i^{\ast} (\alpha) + p^{\ast}(\beta) \ (\beta \in \tilde{K}(\sphere^{2n})) \\
&=& r_0 + r_1 \gamma + \cdots + r_{n-1} \gamma^{n-1} + p^{\ast}(\beta)  \ \mbox{where $r_i \in \D{Z}$} \\
&=& r_0 + r_1 \gamma + \cdots + r_{n-1} \gamma^{n-1} + r \gamma^n.
\end{eqnarray*}

As a result, by induction, we conclude that $\tilde{1}, \gamma, \cdots, \gamma^n$ generate $K(\PC^n)$ over $\D{Z}$ as a ring as desired.

\end{proof}

\nocite{*}
\bibliographystyle{alpha}

\end{document}